\documentclass[11pt]{amsart}
\usepackage{amsmath}
\usepackage{amssymb}
\usepackage{latexsym}
\usepackage{verbatim}
\usepackage{graphics}
\usepackage{graphicx}

\usepackage[all]{xy}

\begin{document}

\title{Now that we're together*: Biography of the Chilean Collective of Women Mathematicians.}
 
\author{Mar\'{\i}a Isabel Cortez and Andrea Vera-Gajardo}
\address{Facultad de Matem\'aticas, Pontificia Universidad Cat\'olica de Chile.
 Instituto de Matem\'aticas, Universidad de Valpara\'{\i}so.} \email{maria.cortez@mat.uc.cl,
andreaveragajardo@gmail.com}
\thanks{Both authors were supported by proyecto ANID PIA Anillo SOC180025.}
\thanks{* "Now that we’re together, now that we’re seen" (\textit{Ahora que estamos juntas, ahora que s\'{\i} nos ven}) is a slogan that has been widely used both in the “Ni una menos” (“Not a single woman less”) movement in Argentina, and in the movement born during the so-called Feminist May of 2018 in Chile. }
  \keywords{Chilean women in math organization, Latin American women in science organization, Chilean math community}

\footnotesize{This is submitted to an upcoming book aimed at a general audience. Below is the second version of the draft.}


\begin{abstract} In this article we give an account of the history of the Chilean collective of women
mathematicians. Unlike the AWM, this collective is not a membership organization, but rather a group that was spontaneously armed with common interests, and whose number of members varies between 15 and 5 members, depending on the moment. We begin by describing the context of the mathematical community
in Chile and the process of forming the Collective, together with the first objectives
we set ourselves. Then we continue with an analysis of some reasons that support
the need to create a group formed by women mathematicians and also the choice of
horizontality and autonomy as structural pillars of our organization. On the
other hand, we refer to the main activities that we have carried out and provide an
overview of the women in mathematics organizations that exist in Latin America.
Finally, we conclude by discussing the research project to which some of the
members of the collective are dedicated nowadays in 2020 (in particular the authors of this
article).
\end{abstract}

 \maketitle{}
 

\section*{The context of the mathematical community in Chile.}
\label{sec:1}

According to the International Mathematical Union, the modern development of mathematics in Chile began in 1930 with the arrival of Carlos Grandjot \cite{DeLaPena2014}, a German mathematician who played a key role training the first Chilean mathematicians \cite{Gutierrez2002}. On the other hand, the first local articles published in international journals were written in the 1940s by Roberto Frucht, another German mathematician based in Chile, who worked at the Universidad Técnica Federico Santa María (UTFSM) until the end of his professional career \cite{Mercado2020}. It is worth noting that in 1949 Frucht published "Graphs of degree three with a given abstract group" in the \textit{Canadian Journal of Mathematics}, the same year in which women’s right to vote in presidential elections was approved in Chile.

By 1970, Chilean mathematicians had already graduated from PhD programs abroad. The first was Jaime Michelow, who obtained his PhD from the University of Washington in 1962. Nevertheless, as the study commissioned by the National Commission of Science and Technology (CONICYT) points out  \cite{Harding1970}, by 1970 there was still no "systematic and planned mathematical research, nor high-level researchers to initiate it" \cite[p. 2]{Harding1970}, and there was also "a lack of human resources in terms of quantity, capacity and adequate knowledge, for university teaching and postgraduate studies" \cite[p. 2]{Harding1970}. This reality stood in stark contrast to other Latin American countries, such as Brazil and Argentina. The case of Brazil is emblematic, since by the 1950s it had already trained mathematicians who would later gain international recognition through their guest lectures at the International Congress of Mathematicians \cite{DeLaPena2014}. On the other hand, Argentina had established the Argentine Mathematical Union (UMA) in 1936, which was created with the aim of "promoting the evident progress of mathematical research in Argentina"  \cite{UMA2020}.

In Chile, institutional development for the training of mathematicians began in 1962 with the creation of the first degree in Mathematics at the Universidad de Concepción \cite{Harding1970}. This was followed in 1965 by degrees in Mathematics at both the Universidad Católica and Universidad de Chile, as well as a degree in Mathematical Engineering from the Faculty of Physical and Mathematical Sciences (FCFM) at the latter institution \cite{Harding1970}. In 1967, the Universidad Técnica del Estado (UTE)—now the Universidad de Santiago de Chile (USACH)—created the country's first postgraduate program in mathematics. This corresponded to the Academic Mathematics Degree (LAM) program, which awarded a Master's degree in Mathematics. The first PhD in Mathematics offered in Chile was awarded in 1975 by the Universidad Católica (UC) to Pablo Salzberg.\footnote{Private communication with Irene Mikenberg} The second one was awarded in 1978 by the same program to Irene Mikenberg, who became the first woman to obtain her PhD in mathematics in Chile. By the end of the 1970s, Chile had 47 PhDs in mathematics working in the country, of which five were women: Carmen Cortázar (UC), Irene Mikenberg (UC), Lidia Consiglieri (U. Católica de Valparaíso), María Angélica Astaburuaga (UTE) and Nancy Lacourly (FCFM U. de Chile)  \cite{CONICYT1978}. It is important to highlight Inés Harding, since she played a very relevant role in the mathematical community of the time. Harding read further postgraduates courses leading to a PhD degree at the University of Moscow in the 1960s  \cite{CONICYT1978}, working as a professor \cite{Pedemonte2017} at UTE until 1981.\footnote{Private communication with Rafael Labarca.}  She became one of the pillars of the LAM program \cite{Calderon2003}, managing the program and recruiting students, while also co-authoring  groundbreaking textbooks from those years, such as "Elementos de Computación" (Computing Elements) in 1973 \cite{Michellow1973}. She passed away in Paris in 1996.

One of the many consequences of the 1973 coup d'état was the loss of many students and academics who were forced to leave the country \cite{SOMACHI2014}, thus delaying the consolidation of the local mathematical community. This situation slowly began to turn around in the early 1980s as several mathematicians who had obtained their PhDs abroad started returning to Chile. Founded in 1975, the Chilean Mathematical Society (SOMACHI) was instrumental in drawing these mathematicians back home as well as promoting mathematical activity at the time \cite{SOMACHI2014}.

A milestone that shaped the way science is done in Chile was the creation of the National Fund for Scientific and Technological Development (FONDECYT) in 1981, which finances basic scientific research in Chile through "open competition for public grants" \cite{FONDECYT2020}. Even today, being awarded a FONDECYT project is considered a decisive milestone in a researcher's career. In most cases, opportunities for national and international collaboration, the possibility to train students, and even the academic workload allocated to each person depend on a research project being awarded one of these grants.

Between 1982 and 1990, FONDECYT funded 144 math projects,\footnote{FONDECYT provided us with data for all the FONDECYT mathematics projects between 1983 and 2019.} of which 13 were awarded to women. These women included Carmen Cortázar, Marta García-Huidobro and Rubí Rodríguez—recognized scholars who remain scientifically active to this day— and Cecilia Yarur, who has recently retired after leading a successful career as a researcher and teacher.

After the return to democracy in 1990, the development of Chilean mathematics started to gain momentum. Not only did the number of researchers and publications grow, but also the number of international cooperation schemes \cite{SOMACHI2014}. The number of FONDECYT projects in mathematics awarded in the 1990s grew by more than 80\% compared to the previous decade, although the ratio of women among the researchers behind these projects remained largely unchanged.

There are currently more than 350 mathematics researchers working in Chile  \cite{SOMACHI2014}. According to data collected from 11 mathematics departments in the country, approximately 20\% of the scholars who hold a PhD in mathematics or related fields are women \cite{Cortez2016}. Nevertheless, a review of the statistics of FONDECYT mathematics projects shows that grants awarded to women between 2011 and 2020 \footnote{The numbers about FONDECYT projects since 2001 include both Regular and Initiation calls. Regular calls are oriented to senior researchers whereas Initiation calls are for junior researchers.}  did not surpass 15\%—an even lower rate than the 2001-2010 period, at 17\%. This suggests that women in our community are even more underrepresented in leading positions, \footnote{Being awarded a FONDECYT grant is considered a sort of symbolic credential of prestige for the beneficiary researcher.}and this situation does not seem to have improved over time. To further emphasize this situation, it is worth mentioning that in the "FONDECYT regular" \footnote{This is the fund oriented to senior researchers.} call for 2020 there were no  projects awarded for women in the area of mathematics, something not seen since 1990.   At the school level, the picture is no more promising, as can be seen from the gender gap in performance in mathematics, present since primary school \cite{Radovic18}.  

There seems to be a consensus within our community on the need to increase female participation in the discipline. What is less clear, however, is how to do this; and even more so, whether that is indeed the problem or rather a symptom of something that goes much deeper. Clémence Perronnet clarifies this problem by stating that institutions expect to receive "individualizing and psychologizing" explanations for inequalities in science, which leads to hiding and concealing "the collective and structural dimension" of these inequalities \cite{Perronnet2020}. In other words, the lack of women in mathematics is not a problem of motivation of individuals (women in this case), but rather the consequence of practices that are usually incompatible with the way in which those individuals are socialized.  


\medskip

In the following lines we present the process of forming the Chilean Collective of Women Mathematicians and we refer to the main activities that we have carried out.   We also describe  other  mathematics organizations  in Latin America that support girls' and women's participations in mathematics.

\section*{How the Collective began.}
\label{sec:2}

  After one of the roundtable discussions at the International Congress of Women in Mathematics 2014 in Seoul, several colleagues from Chile who attended met to discuss the issues presented at the conference. The general consensus that day was that our experience had been very different from that of the speakers, in the sense that our parents had not prevented us from studying mathematics, nor had we endured explicit complaints from our husbands or partners about the profession being unladylike. Yet we were familiar with the experience of being the "only woman" in our workplace, and we perceived that treading through a heavily male-dominated space had consequences that were not easy to identify and describe, given how subtle they generally were. Trying to understand this phenomenon triggered our need to meet as a group for further discussion; this group would later become the Chilean Collective of Women Mathematicians.

Our discussion group was structured around weekly meetings in the Council Room of the Mathematics Department in the Science Faculty at Universidad de Chile. Javiera Barrera, M. Isabel Cortez, Leslie Jiménez, Monica Musso, Adriana Piazza, Anita Rojas, Mariel Sáez and Andrea Vera participated in these first meetings. This space was used to share our experiences, as well as the literature that we came across. We think it was very liberating for all of us to realize that we shared the same experiences as many other women, some of which had even been documented. Shortly after we created our discussion group, we received an invitation from USACH students to speak about women in mathematics at an activity organized for the anniversary of the Mathematical Engineering degree. This was a great opportunity to present the issues we discussed every week to a wider audience. Instead of preparing a simple slideshow with biographies of female mathematicians, we decided that it would be much more interesting to shine a light on the situation of women and mathematics in Chile in quantitative terms, with the aim of starting a conversation. The activity was titled "Panel Discussion: Women and Mathematics" (Conversatorio: mujeres y matemáticas), and was held on November 4, 2014, in a nearly full auditorium with capacity for 100 people. The audience included mathematical engineering students (mostly male), some mathematical colleagues, and several social science researchers specializing in gender issues. After our brief presentation, the discussion with the audience lasted almost two hours. Students' stories related to their experience and difficulties of studying in a place where there are mostly men. Some colleagues were uncomfortable with questions about the Math Olympiad that arose during the conversation, while colleagues from the social sciences summarized and organized the discussion within certain theoretical frameworks. What took place on that day made us realize that there was a great need to address certain issues. We agreed that meeting once a week to discuss these matters was not enough. That is when we decided to come together under the umbrella of the Chilean Collective of Women Mathematicians, as a way of making gender issues in mathematics and our future efforts visible. Having the experience of other women's mathematical associations in different countries as a reference—for example, 
the Association for Women in Mathematics (AWM)—was crucial to creating the Collective. 

It is important to mention that the collective is not the only initiative carried out by  women in mathematics in Chile. Some of these examples are MATEA\footnote{https://matea.cl/hypatia/}  and Colectivo Resistencia Notheriana\footnote{https://fiestadelaciencia.cl/colectivo-resistencia-noetheriana-universidad-de-chile/}.  The first one focus on the organization of summer maths camps for high school girls, and the second one aims to ensure inclusion in the teaching and dissemination of mathematics, from an interdisciplinary perspective. One of the founders of Resistencia Notheriana is also a member of the Collective. Besides mathematics, it is important to highlight the Red de Investigadoras\footnote{https://redinvestigadoras.cl} (Network of Researchers), which is an association that promotes gender equality in research in all areas of knowledge.  A mathematician who is part of the Collective of Women in Mathematics was its first president.

\section*{Organizational characteristics of the Collective.}
\label{sec:3}

Unlike other organizations, the objectives of our Collective were not entirely clear upon its establishment. In fact, the main reason why we decided to band together was the need to be together. We needed to share our experiences with other women mathematicians and confirm that, despite the differences in our lives and academic careers, the sole condition of being women mathematicians in Chile was a shared experience that united all of us. Every one of us had many questions. In time, we realized the need to share them with other women mathematicians, so we could begin to understand them and seek answers. The decision to form an all-women's collective can be explained through Julieta Kirkwood's words: "Only women have the problem of asking themselves what their condition means. Womanhood is otherness (...)" \cite{Kirkwood1987}.

It is important to say that unlike the AWM, this Collective is not -yet- a membership organization, but rather a group that was spontaneously conformed with common interests, and whose number of members varies between 15 and 5 members, depending on the moment and the particular initiative that is being carried out. 

We were very clear from the start about the way in which we made our decisions and the structure that we wanted for our organization. We all agreed that our organization was to be horizontal; in other words, that its organizational structure should be flat, without hierarchies, resulting in democratic decision making which means most of the time by consensus and in some other times by simple majority. As stated by Pérez-Arrau, Espejo, Mandiola, Ríos and Toro \cite{Perez-Arrau2020}, "groups lacking allocated hierarchical leadership may develop different forms of collective leadership." Much like these authors, our notion of collective leadership is a relationship that is built among us, one that is nourished by taking advantage of our different paths and experiences, and that gives way to collective decision making. On the other hand, since our respective workloads are different and dynamic, our time availability for the Collective also varies. Therefore, we could also talk about a kind of "occasional leadership," in which depending on the time available and what each individual seeks, members may from time to time take a leading role to organize specific activities. Still, there are certainly many initiatives that have been possible thanks to the commitment and dedication of everyone involved. We will talk about this in detail in the next section.

Although our work takes place within academic institutions, we decided from the onset to be totally independent, both organizationally and financially. That is to say, we do not want to depend on any institutions, including the Mathematical Society of Chile. The purpose of this decision was to remain free to discuss whatever we believe to be important, all the while protecting our capacity for independent thought. However, from time to time, especially when organizing large events, we resort to some financing and contingent sponsorship. For example, to organize the II Meeting of Women in Mathematics in Latin America, held in 2018 in Valdivia, Chile, we received funding from various institutions such as: International Mathematical Union's Committee for Women in Mathematics (CWM), Conicyt Chile, our home universities, and some local research grants in mathematics.

Today we can say that these decisions we made for our organization—horizontality and independence—turned out to make our Collective unique compared to similar organizations across Latin America. The establishment of the Collective, however, has been questioned by some colleagues in the local mathematical community. During our first year, we sensed that few colleagues understood the need for and relevance of an organization devoted to issues about gender and mathematics. We even saw considerable resistance from some colleagues to our organization. Moreover, our decision to be horizontal and independent also brought about challenges and greater efforts to achieve credibility within the local community. This is definitely something we have in common with AWM and its history, sharing the same feelings expressed in \cite{Greenwald2015} when referring to the decision made by AWM to be autonomous from AMS and MAA: "the advantages of forming an independent association outweighed the challenges of doing so."

Another principle that we share with AWM is the desire to foster a cultural change within the mathematics community, so that young people who are new to the field will experience a less hostile environment than those who came before them. We also seek to highlight and put a spotlight on the life and work of women throughout Latin America who have devoted themselves to mathematics. In fact, one of our permanent activities is a photo exhibition of women mathematicians from Latin America and Europe, which we will explore in the following section.

\section*{Joys and challenges: some activities organized by the Collective.}

The Maths and Gender Panel Discussion (Conversatorio Matemáticas y Género) held on November 27, 2015 in Pucón, has been one of the most complex activities we have organized so far, mainly because of the high amount of exposure it meant for us. This discussion was held as part of SOMACHI's annual meeting and it was probably the first time that gender issues in mathematics were addressed at such a massive local scientific event. Participants included Javiera Barrera, M. Isabel Cortez, Adriana Piazza, Anita Rojas and Mariel Sáez, in addition to social science researcher Marcela Mandiola, who we invited as a guest speaker after having met her at the first discussion group at USACH. Organizing this event was far from easy. In fact, a dispute about which names to include on the list of sponsors in the promotional poster even led to a rift within the Collective. Yet the activity carried on in spite of it all, generating enormous expectation among the local mathematical community. On the day of the event, during lunchtime, we realized that the congress participants (mostly men) felt a kind of "moral duty" to attend our discussion, and that not doing so was somehow considered socially reprehensible. And indeed, when the event took place, everyone at SOMACHI was there, including foreigners who could hardly speak any Spanish, which was the language in which the discussion took place. The activity consisted of two lectures: in the first one, M. Isabel Cortez provided a quantitative overview of women mathematicians in Chile. The second lecture was given by Marcela Mandiola, dealing with her research paper on gender in academia \cite{Mandiola19}.  In her presentation she explained how the different practices in academia where hierarchically gendered, where research and teaching appear as masculinized and  feminized respectively. Attendees found these conclusions to be very controversial, triggering a discussion that lasted for hours. Some colleagues felt under attack, others kept making jokes in mathematical language, and some others expected to hear concrete solutions. The conversation ended very late, so we all went to have dinner together. The discussion went on all through the evening and until the end of the congress. For us, the greatest achievement of this activity was to have introduced the topic to our colleagues.

In 2017 Mariel Sáez suggested the possibility of bringing to Chile the travelling photographic exhibition "Women of mathematics throughout Europe: a gallery of portraits," created by mathematician Sylvie Paycha and photographer Noel Matoff and presented for the first time in 2016 at the European Congress of Mathematics. As a Collective we approved of the idea but decided to modify the original exhibition to also include portraits of women mathematicians from Latin America. Thus, the Mathematics Portraits Exhibition was born, whose first version took place between March 8 and 28, 2018, at the GAM Cultural Center in Santiago, Chile. This was a major project that involved many people and received funding from different sources. The first stages were managed by CONICYT's EXPLORA program. The interviews were conducted by Paula Arenas and the photographs were taken by Cristian Translaviña, Nicolás Sanhueza and Dario Cuellar.  We chose to portray diversity in mathematics from multiple perspectives: regional (not only mathematics in the capital), institutional (from different universities), age and also in terms of their main occupation, taking into account that we were also interested in showcasing outstanding women who excel not only in their research, but also in terms of teaching or outreach initiatives. This has been the largest activity we have organized for a more general target audience. The exhibition has continued to travel around Chile and Latin America (Uruguay), currently thanks to the funding provided by a research project involving M. Isabel Cortez, Mariel Sáez and Andrea Vera. Thanks to this exhibition Mariel Sáez edited a book called "Retratos de Matemáticas" (Mathematical Portraits)  \cite{Saez2018}, whose printed version will be published soon.

Organizing and running the Second Meeting of Latin American Women in Mathematics has been undoubtedly our most important event, and probably the one we are most proud of. It all started at the meeting dubbed "Women in Mathematics in Latin America: Barriers, Advancements and New Perspectives," held in August 2016 at the Casa Matemática de Oaxaca (BIRS), Mexico. M. Isabel Cortez attended this meeting and suggested Chile as the venue for the next such encounter. The project was kicked off with an application to the IMU's Committee for Women in Mathematics (CWM) 2017 call for proposals, which we fortunately won. Even though the awarded amount was not enough for the event, it still meant a lot to us because it encouraged us to shape our expectations and venture forth to get more funding, form a scientific committee, start inviting the first speakers, etc. The organizing committee was made up of  Javiera Barrera, M. Isabel Cortez, M. Isabel del Río, Adriana Piazza, Salomé Martínez, Mariel Sáez and Andrea Vera. We decided to hold the meeting in the last week of January 2018, which is summer in the Southern Hemisphere and right before the holiday season in Chile. The place: Valdivia, a beautiful city in the South of Chile, bordering Patagonia \footnote{The most important street in Valdivia is called Ram\'on Picarte, in recognition of the father of the first Chilean mathematician (of the same name) having international recognition. \cite{GutierrezGutierrez16}} . While preparing for this meeting, we once again encountered less-than-supportive reactions from some colleagues who did not see the need or the relevance of holding a meeting on Women in Mathematics. It would seem that questioning certain practices within our community is seen as a threat by some.

It was clear to us that in order to deal with the issue of women in the sciences, it is necessary to draw on knowledge and expertise from other areas such as the social sciences and philosophy. Therefore, the event was roughly divided as follows: 50\% of the time was spent on lectures and posters on mathematics (presented by women), 25\% on lectures in the social sciences and/or philosophy and the rest of the time was devoted to a leadership workshop and a roundtable entitled "Mujeres en ambientes masculinizados” (Women in male-driven environments), which was open to the general public. Also, during the whole week and in parallel to the above-mentioned activities, we set up an exhibition called “Iluminación Matemática” (Mathematical Illumination), organized by Mónica Canales.
The fact that we included the participation of social scientists in this meeting was undoubtedly a great decision. These points of view gave new insights to our discussions, which led us to approach the issue of women in science/mathematics from a structural perspective, allowing us to get a taste of interdisciplinary collaboration and its advantages.


 \section*{Women's mathematical organizations in other Latin American countries.}

Talking about the Chilean Collective of Women Mathematicians also involves our colleagues and allies in other Latin American countries, since to some extent we consider them part of our organization. Indeed,  some of the activities where the Collective has been involved were organized  jointly with these organizations. For example, the panel of Gender and Mathematics of the Mathematical Congress of the Americas, held in  2017 in Montreal, was organized  together with the Mexico's organization described below. The  next conference of Women in Mathematics in Latin America (to be held next year in Colombia), is another example  of international collaboration.  

Between 2013 and 2018, several organizations of women mathematicians (and scientists in general) emerged throughout Latin America. Many of these organizations, in addition to the Collective in Chile, constitute what we have named the "Network of Women in Mathematics in Latin America and the Caribbean" (Red de Mujeres Matemáticas en América Latina y El Caribe), which began to take shape during 2016 and was consolidated in 2018 at the International Congress of Mathematicians (ICM) Rio de Janeiro-Brazil. 

From a historical point of view, it should be noted that in 2018 the "WM2: World Meeting for Women in Mathematics" took place as a satellite event of the ICM in Rio de Janeiro, Brazil. This probably contributed to the creation of some of the organizations that we will mention in the following paragraphs. In this same line, we also think it is relevant to mention that in May 2018 the so-called "Feminist May" broke out in Chile. It started out as a university protest movement to complain about the lack of protocols to prevent and sanction cases of sexual harassment and abuse in higher education institutions. The movement gained a lot of traction and finally became a national and citizen-wide protest united under the slogan "Non-Sexist Education." Much has been written about this social movement, see for example \cite{TFS19} and \cite{Zeran18}, or \cite{Dessi18} for an article in english.

Next, we provide a review of the women's mathematical organizations in Latin America that we are currently aware of. We particularly refer to their organizational characteristics. To do this, we relied on the report by Araujo-Pardo \& Vera-Gajardo (2019) and updated the information according to conversations we had with colleagues in Peru, Argentina and Brazil.\footnote{Roxana López Cruz (Peru), Gabriela Ovando (Argentina) and Christina Brech (Brazil)} 
In Argentina there is the Gender Commission of the Argentine Mathematical Union (UMA), founded in 2018. At first it was made up only of women, but now it is a mixed group of 5 people, including a female coordinator. They also have representatives in different parts of the country. The committee sees itself as a coordinating entity rather than a resolving body. It works as a gender advisor to the UMA and proposes its own specific initiatives.

In Colombia there is the Commission on Equity and Gender of the Colombian Mathematical Society \cite{CEGC2017}, created in 2017. It is currently a mixed organization made up of six members, with a female chairperson. One if its objectives is to encourage participation in mathematics by women and minorities in Colombia.

Mexico’s organization is the oldest in Latin America. It is called the Commission on Equity and Gender of the Mexican Mathematical Society  \cite{CEG2020}, founded in 2013. It is a mixed group of three people, one of whom acts as Coordinator. Its purpose is to promote the inclusion of underrepresented groups, particularly women, in the country's mathematical endeavors. 

The case of Brazil is slightly different. Because it is such a large and diverse country, there are different organizations and networks. Some of them include the Gender Commission \cite{CDG2020} created in 2018 as a joint effort between the Brazilian Mathematical Society and the Brazilian Society for Applied and Computational Mathematics. On the other hand, some Brazilian women mathematicians created the platform "Mathematics, a feminine noun" \cite{MSF2020}, which aims to promote actions being developed in different parts of the world, especially in Brazil, in order to increase women's participation in mathematics. We should underscore that these are just some of the many networks that exist in Brazil.

In Peru, the "Group of Peruvian Women in Mathematical Sciences" \cite{MPM2016}  was formed in 2016. This is a network of more than 100 members, which has an organizing committee, a list of representatives in different parts of the country and financially dependent on the Peruvian Society of Applied Mathematics and Computer Science. The main—although not the only—activity performed by this group is to coordinate and organize the Sofia Kovalevskaia Prize in Peru, awarded in many countries by the eponymous Fund.

\section*{The present and the future: paths taken by the Collective’s members.}

In 2020, the members of the Collective continue to carry out some activities together and others on their own, focusing on topics related to mathematics and gender. For example, we continue to regularly display the exhibition "Retratos de Matemáticas" (\textit{Mathematical Portraits}) in different cities and contexts, as well as attending most of the activities in which we are invited to present. We are also in permanent touch with our partners in the "Network of Women in Mathematics in Latin America and the Caribbean."

Javiera Barrera and Adriana Piazza generally participate in a few national activities to spread scientific knowledge, such as the "Feria Ingeniosas." Similarly, the Collective is one of the organizations involved in the global coordination of the "12 de Mayo Celebration of Women in Mathematics initiative honoring Maryam Mirzakhani on her birthday.", represented by Andrea Vera. Javiera Barrera is also part of the Equity and Gender Committee of the Faculty of Engineering and Sciences at her home institution.

Since 2018, M. Isabel Cortez, Mariel Sáez and Andrea Vera have been part of an interdisciplinary research team in the "Anillo Matemáticas y Género" (Mathematics and Gender “Anillo” Grant)  \cite{Anillo2020}. The team is also integrated by Jeanne Hersant, Marcela Mandiola and Tania de Armas, who are all social science researchers. This is a three-year project funded by the Chilean National Agency for Research and Development. The main objective of this research is to understand and analyze the situation of women in mathematics in Chile from a gender perspective, told through the stories and career paths of those who are part of it. To achieve this, we have separated it into three specific goals giving rise to four stages of research. These goals are as follows:
\begin{enumerate}
\item To gain awareness of the scientific field of mathematics in Chile from a gender perspective, as well as its subsequent analysis. 
\item To understand and analyze the subjectivities of female mathematics scholars from their interpretations
and self-interpretations regarding the academic field of mathematics, as well as in relation to their career paths.
\item To identify the common career pattern of women researchers in mathematics, that is, the different stages they must go through to reach the status of researcher; in addition, to analyze the barriers and difficulties that women researchers in mathematics must overcome to become researchers, in relation to internal hierarchies and prestige ladders.
\end{enumerate}

Regarding the future, in these months (November, December 2020) we are preparing a workshop that will have the objective of re-founding the Collective. Specifically, in January 2021 we will hold an open day to women in mathematics in Chile, who wish to be part of the Collective. Additionally, one of our plans is to create, within the Collective, a chapter for students.

 On the other hand, together with AWM and some Mathematician colleagues from Latin America, we are organizing a virtual bilingual (English-Spanish) math outreach program. 

Undoubtedly,  since the formation of the Collective, our careers have been strongly influenced by issues that touch upon gender and mathematics. Without a doubt, one of the lessons learned from belonging to a Collective of Women Mathematicians is the certainty that any change we want to promote requires coordination with others; that is, collaboration.


\section*{Acknowledgement}   We would like to thank Hernán Henríquez, Rafael Labarca and Irene Mikenberg for providing us valuable information regarding the history of mathematics in Chile. We would also like to thank FONDECYT for providing us the complete database concerning FONDECYT grants in mathematics.

\end{document}